%% file: arXiv.tex
\documentclass[a4paper,10pt]{article}
\usepackage{paper-en}
\usepackage{hyperref}

\def\thetitle{Quadratic metric comparisons}
\def\theauthors{Nina Lebedeva, Anton Petrunin, and Vladimir Zolotov}

\hypersetup{colorlinks=true,
citecolor=black,
linkcolor=black,
anchorcolor=black,
filecolor=black,
menucolor=black,
urlcolor=black,
pdftitle={\thetitle},
pdfauthor={\theauthors}
}

\begin{document}

\title{\thetitle}
\author{\theauthors}
\date{}
\maketitle

\input{4-point}

{\sloppy
\def\emph{\textit}
\printbibliography[heading=bibintoc]
\fussy
}
\end{document}

%% file: 4-point.tex
\begin{abstract}
We investigate how quadratic constraints on the six pairwise distances in every four-point array affect the geometry of length spaces.
\end{abstract}

\section{Introduction}\label{par:quadratic-inq}

\paragraph{Quadratic condition.}
Let $\bm{x}\z=(x_1,\ldots,x_n)$ be an $n$-point array in a metric space $X$ and let $a_{i,j}$ be the components of a symmetric $n{\times}n$ matrix.
An inequality of the following type
\[\sum_{i,j}a_{i,j}\cdot|x_i-x_j|_X^2\geqslant 0\]
will be called \emph{quadratic}.
A system of quadratic inequalities will be called a \emph{quadratic condition}.

We will be interested in length spaces $X$ such that a given quadratic condition holds for any $n$-point array in $X$.
The following statement says that the Toponogov theorem has no simple relatives.
It will be proved in Section~\ref{par:globalization}.

\begin{thm}{Main theorem}
Suppose that a quadratic condition for 4-point arrays satisfies the globalization property;
that is, if this condition holds locally (in a neighborhood of any point) in a length space $X$, then it holds in $X$.
Then either the condition is trivial, or it describes Alexandrov spaces with nonnegative curvature.
\end{thm}

\paragraph{Auxiliary results.}
Let us state two auxiliary results that are of independent interest.
Their formulations use a special case of quadratic inequalities
with $a_{i,j}=-\lambda_i\cdot\lambda_j$ for a real array $(\lambda_1,\ldots, \lambda_n)$ such that $\lambda_1+\ldots+\lambda_n=0$.
These are so-called inequalities of negative type; they will be discussed in Section~\ref{par:rank-one}.

In Section~\ref{Four-point arrays}, we prove the following version of a theorem by Abraham Wald \cite[§ 7]{wald} describing the metrics of all possible 4-point arrays in Alexandrov spaces with nonnegative and nonpositive curvature.

\begin{thm}{Proposition}
The following three conditions are equivalent:
\begin{enumerate}
\item A 4-point metric space $X$ is isometric to a subset of a length space with nonnegative (nonpositive) curvature in the sense of Alexandrov.
\item $X$ is isometric to a subset of the product $\RR^3\times r\cdot \SSS^1$ for some $r>0$ (respectively, of $\RR^3\times Y$, where $Y$ denotes the tripod; that is, three half-lines with a common base point).
\item All inequalities of negative type with $\lambda$-arrays such that $\lambda_1\cdot\lambda_2\cdot\lambda_3\cdot\lambda_4<0$ (respectively, with $\lambda_1\cdot\lambda_2\cdot\lambda_3\cdot\lambda_4>0$) hold in $X$.
\end{enumerate}
\end{thm}

The next proposition says that a single nondegenerate inequality of negative type defines Alexandrov spaces with either nonnegative or nonpositive curvature; it will be proved in Section~\ref{Alexandrov's comparison}.

\begin{thm}{Proposition}
Let us consider all length spaces that satisfy the quadratic condition given by an inequality of negative type with $\lambda$-array $(\lambda_1,\lambda_2,\lambda_3,\lambda_4)$.
\begin{enumerate}
\item If $\lambda_1\cdot\lambda_2\cdot\lambda_3\cdot\lambda_4<0$, then it describes all Alexandrov spaces with nonnegative curvature.
\item If $\lambda_1\cdot\lambda_2\cdot\lambda_3\cdot\lambda_4>0$, then it describes all Alexandrov spaces with nonpositive curvature.
\item If $\lambda_1\cdot\lambda_2\cdot\lambda_3\cdot\lambda_4=0$, then it gives no restrictions; it describes all length spaces.
\end{enumerate}

\end{thm}

This statement shows that a relatively weak form of 4-point comparison forces a much stronger 4-point comparison (once we assume that the metric is intrinsic).
The proof is a straightforward combination of an argument by Takashi Sato \cite{sato} and a variation due to two of the authors \cite{lebedeva-petrunin-2010};
these papers consider particular cases of such inequalities for the $\lambda$-arrays $(1,1,-1,-1)$ and $(1,1,1,-3)$.

Sections \ref{Associated form} and \ref{par:rank-one} introduce necessary definitions.
Here we define the associated quadratic form for a point array and the inequalities of negative type.
We also prove several basic statements that are valid for all $n$-point arrays.
Section \ref{Auxiliary statements} provides a technical statement for the proof of the main theorem.

\section{Associated form}\label{Associated form}

Choose a point array $\bm{x}=(x_1,\ldots,x_n)$ in a metric space.
Let $V_n\z=\RR^{n-1}$ be the Euclidean space with a choice of regular simplex $\triangle$, whose vertices are $y_1,\ldots,y_n$.
Consider the quadratic form $\rho_{\bm{x}}$ on $V_n$ defined by the equalities
\[\rho_{\bm{x}}(y_i-y_j)=|x_i-x_j|^2_X\]
for all indices $i$ and $j$.

The quadratic form $\rho_{\bm{x}}$ will be called the \emph{associated form} of the array $\bm{x}$;
it is uniquely defined and remembers all distances $|x_i-x_j|_X$
(we assume that the simplex $\triangle$ in $V_n$ is known).
We are allowed to not distinguish between $\rho_{\bm{x}}$ and the corresponding semimetric on $\bm{x}=(x_1,\ldots,x_n)$.

The Euclidean structure on $V_n$ identifies it with $V_n^*$.
The space of quadratic forms on $V_n$ will be denoted by $W_n$;
it is the symmetric square of $V_n=V_n^*$, and can be written as $W_n=S^2(V_n^*)=S^2(V_n)$.

Any quadratic inequality described above defines a linear inequality on $W_n$, so it can be written as $\langle\omega,\rho_{\bm{x}}\rangle\geqslant 0$ for a fixed $\omega\in W_n$.
Therefore, any quadratic condition defines a closed convex cone in $W_n$, say $K$;
that is, $K$ is a nonempty closed set such that if $v,w\in K$, then $a\cdot v+b\cdot w\in K$ for any $a,b\geqslant0$.
Denote by $\mathcal{M}_K$ the class of all length spaces $X$ such that
$\rho_{\bm{x}}\in K$ for any $n$-point array $\bm{x}\z=(x_1,\ldots,x_n)$ in~$X$.
The class $\mathcal{M}_K$ respects distance-preserving embeddings; that is,
if there is a distance-preserving embedding $X\to Y$ between length spaces and $Y \in  \mathcal{M}_K$, then $X\in \mathcal{M}_K$.

The following two observations were made by Alexandr Andoni, Assaf Naor, and Ofer Neiman \cite[1.4.1]{ANN}.

Since $K$ is a convex cone,
\[X,\ Y\in  \mathcal{M}_K
\qquad\Longrightarrow\qquad
X\times Y\in\mathcal{M}_K
\quad\text{and}\quad
a\cdot X\in\mathcal{M}_K
\]
for any $a\geqslant 0$;
here $X\times Y$ denotes the $\ell_2$-product of metric spaces, and
$a\cdot X$ denotes the rescaled copy of $X$ with scaling factor $a$.
Moreover, all forms $\rho_{\bm{x}}$ for point arrays $\bm{x}$ in $\mathcal{M}_K$-spaces form a convex cone.
Let us denote this new cone by $K'$.

Evidently, $K'\subset K$, and this inclusion might be strict.
The first reason comes from the triangle inequality, which is equivalent to $\rho_{\bm{x}}(w)\geqslant 0$ for any vector $w$ lying in a 2-face of $\triangle$.
These inequalities must hold for any form in $K'$.
Furthermore, if $\rho_{\bm{x}}\in K'$, then $\rho_{\hat{\bm{x}}}\in K'$ for any $n$-point array $\hat{\bm{x}}$ chosen from the points of $\bm{x}$; in particular, $\hat{\bm{x}}$ might be a permutation of points in $\bm{x}$.

These two conditions hold for all metrics (not necessarily intrinsic),
and one cannot get more for general metric spaces.
The next observation already uses the assumption that the metric is intrinsic.

\begin{thm}{Observation}
If $K$ is a closed convex cone in $W_n$, then $K'$ is closed.
Moreover, $\mathcal{M}_K$ is closed under ultralimits.
\end{thm}

(For \emph{ultralimits} and \emph{ultracompletions} of metric spaces and all related topics, see, for example, \cite{petrunin2023}.)

\parit{Proof.}
The last statement is evident, and it implies the first statement.

Indeed, for any sequence of spaces $X_n$ in $\mathcal{M}_K$, its ultralimit $X_\omega$ also belongs to $\mathcal{M}_K$.
Therefore, given a sequence of point arrays $\bm{x}_n$ in $X_n$,
its ultralimit $\bm{x}_\omega$ in $X_\omega$ has the limiting distances between corresponding points.
Hence, the result.
\qeds

\begin{thm}{Proposition}\label{prop:Associated form}
Let $K$ be a closed convex cone in $W_n$.
If $\mathcal{M}_K$ is not trivial (that is, it contains a space with at least two points), then $\mathcal{M}_K$ contains all Euclidean spaces.
\end{thm}

\parit{Proof.}
Let $X$ be a space in $\mathcal{M}_K$ with two distinct points.
Consider the ultracompletion $X^\omega$ of $X$;
the space $X$ will be considered as a subset of $X^\omega$.
Since $K$ is closed, the observation implies that $X^\omega\in \mathcal{M}_K$.

Since $X$ is a length space, $X^\omega$ has to be geodesic.
Since $X^\omega$
contains a pair of distinct points, it must contain a nontrivial geodesic.
It follows that $\mathcal{M}_K$ contains a line segment.
By rescaling the segment and passing to the ultralimit, we get $\RR\in \mathcal{M}_K$;
taking products of real lines then yields the result.
\qeds

Let us denote by $Q$ the cone of nonnegative quadratic forms in $W_n$.
Its dual cone $Q^*\subset W_n^*=S^2(V_n)$ is generated by tensor squares of vectors in~$V_n$.

\begin{thm}{Corollary}
Let $K$ be a closed convex cone in $W_n$ such that $K=K'$.
Then either $K\supset Q$ or $K$ is trivial; that is, $K=\{0\}$.
\end{thm}

\parit{Proof.}
Since $\RR\in \mathcal{M}_K$, we get that $\sigma^2\in K$ for any linear function $\sigma\colon V_n\to\RR$.
By the spectral theorem, any form in $Q$ can be written as a sum of squares of linear functions, hence the result.
\qeds

\section{Rank-one inequalities}\label{par:rank-one}
Recall that $\rho_{\bm{x}}$ denotes the associated quadratic form for a given point array $\bm{x}=(x_1,\ldots,x_n)$.

\begin{thm}{Observation}\label{obs:rank-one}
A point array $\bm{x}=(x_1,\ldots,x_n)$ is isometric to an array in a Euclidean space if and only if $\rho_{\bm{x}}(v)\geqslant 0$ for any vector $v$.
\end{thm}

An inequality of type $\rho_{\bm{x}}(v)\geqslant 0$ for a fixed vector $v\in V_n$ will be called a \emph{rank-one} inequality.
It belongs to the class of quadratic inequalities.  
In the notation of Section~\ref{par:quadratic-inq}, it means that $a_{i,j}=-\lambda_i\cdot\lambda_j$ for a real array $(\lambda_1,\ldots, \lambda_n)$ such that
$\lambda_1+\ldots+\lambda_n=0$.
Such inequalities are also known as inequalities of \emph{negative type} \cite{deza-lauren}.
If the array $(\lambda_1,\ldots, \lambda_n)$ contains $i$ positive and $j$ negative numbers,
then we say that this is an inequality of \emph{negative type} $(i,j)$.
Since changing the signs of all $\lambda_i$ does not change the inequality, we can always assume that $i\geqslant j$.

\section{Four-point arrays}\label{Four-point arrays}

For $4$-point arrays, we have two interesting types of rank-one inequalities: negative type $(2,2)$ and $(3,1)$.
The type $(1,1)$ is trivial, and
the type $(2,1)$ follows from the triangle inequality.
In fact, all the triangle inequalities (there are 12 such inequalities for 4 points) are equivalent to all inequalities of negative type $(2,1)$ (there are infinitely many of them).

Consider a rank-one inequality $\rho_{\bm{x}}(v)\geqslant 0$; we may assume that $v$ is a unit vector.
Since the sign of $v$ does not change the inequality, we may assume that it lies in a closed hemisphere bounded by an equator in the direction of one of the facets of $\triangle$.
These equators divide the hemisphere into 4 triangles and 3 quadrangles.
The inequality $\rho_{\bm{x}}(v)\geqslant 0$ has negative type $(2,2)$ or $(3,1)$
if and only if $v$ lies in the interior of a quadrangle or a triangle, respectively.
Equivalently, this inequality is of negative type $(2,2)$ if $v$ points from one edge of the tetrahedron $\triangle$ to the opposite edge, and of type $(3,1)$ if $v$ points from a vertex to the opposite facet (up to the sign of $v$).

\begin{wrapfigure}{o}{36mm}
\centering
\vskip-3mm
\includegraphics{mppics/pic-20}
\vskip-0mm
\end{wrapfigure}

If a vector $v$ is parallel to a facet of $\triangle$, then $\rho_{\bm{x}}(v)\geqslant 0$ is an inequality of negative type $(2,1)$, which follows from the triangle inequality.
The picture shows the hemisphere.
The labels on the edges indicate which triangle inequality becomes an equality when $\rho_{\bm{x}}$ vanishes at a vector on that edge;
for example, the label $123$ means that
\[|x_1-x_2|+|x_2-x_3|=|x_1-x_3|.\]
If $\rho_{\bm{x}}$ vanishes on the intersection of equators, then two points in the array have to coincide;
the label shows which pair.
For example, if it is marked by $12$, then $x_1=x_2$.

\begin{thm}{Proposition}\label{prop:Four-point arrays}
Let $X$ be a 4-point metric space.

$X$ satisfies all inequalities of negative type $(3, 1)$ if and only if it admits an isometric embedding into the product $r\cdot \mathbb{S}^1\times\RR^3$ for some $r>0$.

$X$ satisfies all inequalities of negative type $(2, 2)$ if and only if it admits an isometric embedding into the product $Y\times\RR^3$, where $Y$ denotes the \emph{tripod};
that is, three half-lines with a common base point.
\end{thm}

Inequalities of negative type $(3, 1)$ hold in Alexandrov spaces with nonnegative curvature; this follows from the so-called Lang--Schroeder--Sturm inequality \cite{lang-schroeder, sturm}.
Similarly, inequalities of negative type $(2, 2)$ hold in Alexandrov spaces with nonpositive curvature.
This follows easily from the (2+2)-point comparison \cite[9.5]{AKP-2024}.
Since the tripod $Y$ has nonpositive curvature, and $\mathbb{S}^1$ has nonnegative curvature in the sense of Alexandrov, we get the following corollary, which also follows from the result of Abraham Wald \cite[§ 7]{wald}.

\begin{thm}{Corollary}\label{cor:Four-point arrays}
A 4-point metric space $X$ is isometric to a subset of a length space with nonnegative (nonpositive) curvature in the sense of Alexandrov if and only if all inequalities of negative type $(3, 1)$ (respectively, type $(2, 2)$) hold in $X$.
\end{thm}

The five-point versions of this corollary have been proved by two of the authors \cite{lebedeva-petrunin-2024} and Tetsu Toyoda \cite{toyoda,lebedeva-petrunin2021}, respectively.

\parit{Proof of \ref{prop:Four-point arrays}.}
Let us enumerate the points in $X$ and let $\rho$ be the associated form on~$V_4$.
By \ref{obs:rank-one} we can exclude the case $\rho\geqslant 0$.

Choose a minimal form $\tilde\rho\leqslant \rho$ such that all $(2,1)$-inequalities hold for $\tilde\rho$;
here $\tilde\rho\leqslant \rho$ means that $\tilde\rho(v)\leqslant \rho(v)$ for any vector $v$.
Consider the (semi)metric on $X$ with the associated form $\tilde\rho$;
denote the corresponding metric space by $\tilde X$.

By \ref{obs:rank-one}, $X$ isometrically embeds into $\tilde X\times \RR^3$.
Hence, it is sufficient to show that $\tilde X$ embeds into $r\cdot \mathbb{S}^1$ for some $r>0$, or, respectively, into $Y$.

Let $N\subset \mathbb{S}^2$ be the set on which $\tilde\rho$ is negative and let $\bar N$ be its closure.
Since $\tilde\rho\leqslant \rho$ and $\rho\ngeqslant0$, we have that $N\ne \emptyset$.

Since all $(2,1)$-inequalities hold for $\tilde\rho$, the form $\tilde\rho$
must be nonnegative on 4 equators $e_1,e_2,e_3,e_4$ in the directions of facets of $\triangle$.
Since $\tilde\rho$ is minimal, $\bar N$ has to touch the union of the equators in at least three directions (up to sign).
If not, then there is a linear function, say $\sigma$, that vanishes at all common points of $\bar N$ and the union of the equators.
In this case, consider the form $\tilde\rho-\varepsilon\cdot \sigma^2$ for small $\varepsilon>0$;
note that all $(2,1)$-inequalities still hold for this form.
Therefore, $\tilde\rho$ is not minimal --- a contradiction.

It means that $\bar N$ lies in a quadrangle or a triangle (up to sign) and touches its sides at three points or more.

In the case of a triangle, $\tilde \rho$ and $\rho$ satisfy all inequalities of negative type $(2,2)$, and $\bar N$ has to touch all sides of the triangle.
According to the diagram above, after relabeling, we can assume that
\begin{align*}
|x_1-x_2|&=|x_1-x_4|+|x_4-x_2|,
\\
|x_2-x_3|&=|x_2-x_4|+|x_4-x_3|,
\\
|x_3-x_1|&=|x_3-x_4|+|x_4-x_1|.
\end{align*}
In this case, the array can be embedded into the tripod $Y$.

In the case of a quadrangle, $\tilde \rho$ and $\rho$ satisfy all inequalities of negative type $(3,1)$, and $\bar N$ touches at least three sides of the quadrangle, but might touch all four.
Look at the diagram and convince yourself that after relabeling, we may assume that
\begin{align*}
|x_1-x_4|&=|x_1-x_2|+|x_2-x_4|=|x_1-x_3|+|x_3-x_4|,
\\
|x_2-x_3|&=|x_2-x_4|+|x_4-x_3|.
\end{align*}
If it touches all sides, we also have $|x_2-x_3|=|x_2-x_1|+|x_1-x_3|$.
In any case, the array can be embedded into $r\cdot \mathbb{S}^1$, where $r=|x_1-x_4|/\pi$;
so the points $x_1$ and $x_4$ become antipodal in $r\cdot \mathbb{S}^1$.
\qeds

The picture shows the possible positions of the set $N$.
\begin{figure}[h!]
\centering
\vskip-0mm
\includegraphics{mppics/pic-30}
\vskip-0mm
\end{figure}
Below it, we provide a diagram following the convention from \cite{lebedeva-petrunin-2024};
if three points, say $x_1$, $x_2$, and $x_3$, appear in that order on a smooth line, then $|x_1-x_2|+|x_2-x_3|=|x_1-x_3|$.

\section{Alexandrov's comparison}\label{Alexandrov's comparison}

\begin{thm}{Proposition}\label{prop:Alexandrov's comparison}
Suppose $K\subset W_4$ is defined by a single rank-one inequality on $4$-point arrays.
\begin{enumerate}
\item If the inequality is of negative type $(2,2)$, then $\mathcal{M}_K$ consists of all length spaces with nonpositive curvature in the sense of Alexandrov.
\item \label{prop:Alexandrov's comparison:(3,1)} If the inequality is of negative type $(3,1)$, then $\mathcal{M}_K$ consists of all length spaces with nonnegative curvature in the sense of Alexandrov.
\item\label{prop:Alexandrov's comparison:all} In the remaining cases, $\mathcal{M}_K$ consists of all length spaces.
\end{enumerate}

\end{thm}

This statement and Corollary~\ref{cor:Four-point arrays} imply that a single inequality of negative type $(2,2)$ or $(3,1)$ on a length space implies \emph{all} inequalities of the same type.

\parit{Proof.}
Our inequality can be written as 
\[\sum_{i,j}\lambda_i\cdot\lambda_j\cdot|x_i-x_j|_X^2\leqslant 0,
\eqlbl{eq:lambda}
\]
where $\lambda_1+\lambda_2+\lambda_3+\lambda_4=0$.
If $\lambda_i=0$ for some $i$,
then the inequality follows from the triangle inequality.
In this case, $\mathcal{M}_K$ consists of all length spaces, which is case \ref{prop:Alexandrov's comparison:all}.

It remains to consider the inequalities of negative type $(2,2)$ or $(3,1)$.
In these cases, we can assume that our $\lambda$-array is
\[(\alpha\cdot (1-\beta),\  (1-\alpha)\cdot(1-\beta),\  \beta,\ -1)\] 
for some $\alpha,\beta$ such that $0< \beta< 1$;
in the case of the $(3,1)$-inequality, we have $0<\alpha<1$, and
in the case of the $(2,2)$-inequality, we have $1<\alpha$.

Consider the 4-point array $\bm{x}=(x_1,x_2,x_3,x_4)$  in the plane such that 
\[x_4=\alpha\cdot (1-\beta)\cdot x_1+(1-\alpha)\cdot(1-\beta)\cdot x_2+\beta\cdot x_3.\]
We can assume that the array matches one of the configurations in the pictures below,
so for the $(3,1)$-inequality, the point $x_4$ lies inside the triangle $x_1x_2x_3$,
and for the $(2,2)$-inequality, the segment $[x_1x_3]$ intersects $[x_2x_4]$.

Note that for this array we get equality in \ref{eq:lambda}.
Moreover, there is a bijection between plane $\bm{x}$-arrays up to affine transformation and $\lambda$-arrays without zeros up to multiplication by a nonzero coefficient;
as before, we assume that the $\bm{x}$-arrays are in general position --- no three of their points lie on one line.
In other words, we can describe our inequality by a 4-point array in general position up to affine transformation.
If the points of the array lie at the vertices of a convex quadrangle,
then it corresponds to an inequality of type $(2,2)$.
If one of the points lies inside the triangle formed by the remaining points, then it corresponds to an inequality of type $(3,1)$.

\begin{figure}[ht!]
\vskip-0mm
\centering
\includegraphics{mppics/pic-10}
\vskip0mm
\end{figure}

Consider the affine transformation that sends $x_1\mapsto x_1$, $x_2\mapsto x_2$, and $x_3\z\mapsto x_4$;
suppose $x_4\mapsto x_5\mapsto x_6\mapsto\ldots$, so
\begin{align*}
x_{3+k}&=\alpha\cdot (1-\beta)\cdot x_1+(1-\alpha)\cdot(1-\beta)\cdot x_2+\beta\cdot x_{2+k}=
\\
&=\alpha\cdot (1-\beta^k)\cdot x_1+(1-\alpha)\cdot(1-\beta^k)\cdot x_2+\beta^k\cdot x_3
\end{align*}
for $k\geqslant 1$.
Note that $x_4,\ldots, x_{2+k}$ lie between $x_3$ and $x_{3+k}$ in the same order;
in particular,
\[|x_3-x_4|+\ldots+|x_{2+k}-x_{3+k}|=|x_3-x_{3+k}|.\]

\begin{thm}{Claim}\label{clm:1=>2}
Suppose that inequality \ref{eq:lambda} for $\lambda$-array
\[(\alpha\cdot (1-\beta),(1-\alpha)\cdot(1\z-\beta), \beta,-1)\] holds for any 4-point array in $X$.
Then so does the inequality for the $\lambda$-array
\[(\alpha\cdot (1-\beta^k), (1-\alpha)\cdot(1-\beta^k), \beta^k,-1)\] for any integer $k\geqslant 1$.
\end{thm}

These two inequalities correspond to the arrays $x_1$, $x_2$, $x_3$, $x_4$ and $x_1$, $x_2$, $x_3$, $x_{3+k}$;
let us denote by $L(p, q, r, s)$ and $L'(p, q, r, s)$ the left-hand sides in \ref{eq:lambda} for these two inequalities evaluated
at $p, q, r, s$.
We know that $L\leqslant 0 $ holds for any 4-point array in $X$ and we need to show the same for $L'\leqslant 0$.

Passing to the ultracompletion, we may assume that $X$ is geodesic.
Choose 4 points $x_1,x_2,x_3,x_{3+k}\in X$ and let $x_4,\ldots,x_{2+k}$ be the points on a geodesic $[x_3x_{3+k}]$ that divide in the same proportions as in our plane configuration;
that is
\[|x_3-x_4|\,:\,|x_4-x_5|\,:\,\ldots\,:\,|x_{2+k}-x_{3+k}|=\beta:\beta^2:\,\ldots\,:\beta^k.\]
If we sum up the first inequality for arrays $(x_1,x_2,x_{2+i},x_{3+i})$ with the appropriate coefficients, then we get the second inequality for $x_1,x_2,x_3,x_{3+k}$.

Namely, $x_{3+i}$ divides $[x_{2+i}x_{4+i}]$ in the ratio $1:\beta$.
Therefore,
\[\beta^k\cdot|x_3-x_{3+k}|^2 =\tfrac{1-\beta^k}{1-\beta}\cdot(\beta^k\cdot |x_3-x_4|^2+\ldots+\beta\cdot|x_{2+k}-x_{3+k}|^2).\]
Hence we get the following telescopic sum
\[
\begin{aligned}
L'(x_1,x_2,x_3,x_{3+k})
=
\tfrac{1-\beta^k}{1-\beta}\cdot\biggl(
&\beta^{k-1}\cdot L(x_1,x_2,x_3,x_4)+\ldots
\\
&\ldots+  \beta^{0}\cdot L(x_1,x_2,x_{2+k},x_{3+k})\biggr).
\end{aligned}
\eqlbl{eq:L'=L+L}
\]
In particular, $L \leqslant 0$ implies $L' \leqslant 0$ as required.

Applying the claim several times, we get all the inequalities with $\lambda$-arrays
\[(\alpha\cdot (1-\gamma),\  (1-\alpha)\cdot(1-\gamma),\ \gamma,\ -1),\]
where $\gamma=\beta^k$ for an integer $k\geqslant 1$;
in particular, we get the following inequality
\[
\begin{aligned}
\alpha\cdot (1-\alpha)\cdot(1-\gamma)^2\cdot|x_1-x_2|^2 - \gamma\cdot |x_3-x_4|^2 &+
\\
+(1-\alpha)\cdot(1-\gamma)\cdot\gamma\cdot|x_2-x_3|^2-\alpha\cdot (1-\gamma)\cdot |x_1-x_4|^2&+
\\
+\alpha\cdot(1-\gamma)\cdot\gamma\cdot|x_1-x_3|^2-(1-\alpha)\cdot(1-\gamma)\cdot |x_2-x_4|^2&\leqslant 0
\end{aligned}
\eqlbl{eq:lambda-inq}
\]
for arbitrarily small $\gamma>0$.

Choose $X\in \mathcal{M}_K$.
Let us apply \ref{eq:lambda-inq}
\begin{figure}[ht!]
\vskip-0mm
\centering
\includegraphics{mppics/pic-15}
\vskip0mm
\end{figure}
to a quadruple $x_1,x_2,x_3,x_4\in X$ such that $x_1$, $x_2$, and $x_4$ lie on one geodesic, and we have equality in the $(2,1)$-inequality with $\lambda$-array
\[(\alpha,\  (1-\alpha),\ 0,\ -1);\]
that is,
\[\alpha\cdot (1-\alpha)\cdot|x_1-x_2|^2-\alpha\cdot |x_1-x_4|^2-(1-\alpha)\cdot |x_2-x_4|^2=0.\eqlbl{eq:trig-inq}\]
Passing to the limit as $\gamma\to 0$ in the inequality $\tfrac1\gamma\cdot$\ref{eq:lambda-inq}$\,-\,\,\tfrac{1-\gamma}\gamma\cdot$\ref{eq:trig-inq}, we get
\[
\begin{aligned}
\alpha\cdot|x_1-x_3|^2+(1-\alpha)\cdot|x_2-x_3|^2-
\alpha\cdot (1-\alpha)\cdot|x_1-x_2|^2 \leqslant |x_3-x_4|^2.\end{aligned}
\eqlbl{eq:CBB-CBA}
\]

Note that if this inequality holds for all $\alpha\in (0,1)$ (or $\alpha\in (1,\infty)$) then we get
a point-on-side comparison for nonnegative (respectively nonpositive) curvature in the sense of Alexandrov \cite[8.14 and 9.14]{AKP-2024}.
Indeed, for Euclidean space, equality in \ref{eq:CBB-CBA} holds for any $\alpha$.
Therefore, if $\alpha<1$, we get $|x_3-x_4|\geqslant |\tilde x_3-\tilde x_4|$, where $\tilde x_4$ divides the side $[\tilde x_1\tilde x_2]$ of the model triangle $[\tilde x_1\tilde x_2\tilde x_3]=\tilde\triangle(x_1x_2x_3)$ in the same ratio $(1-\alpha):\alpha$.
Similarly, if $\alpha>1$, we get $|x_3-x_1|\leqslant |\tilde x_3-\tilde x_1|$, where $\tilde x_1$ divides the side $[\tilde x_4\tilde x_2]$ of the model triangle $[\tilde x_4\tilde x_2\tilde x_3]=\tilde\triangle(x_4x_2x_3)$ in the same ratio $(\alpha-1):1$.

So far, we have obtained an inequality only for one value of $\alpha$;
however, by applying iteration, we can derive inequalities for the remaining values $\alpha\in (0,1)$ (respectively, $\alpha\in (1,\infty)$).
More precisely, assume this inequality holds for some $\alpha\in (0,1)$.
Then we can change $\alpha$ to $(1-\alpha)$, $\alpha^2$,  $\alpha\cdot (1-\alpha)$, $(1-\alpha^2)$, and so on.
In other words, if $\mathcal{A}\subset(0,1)$ is the set of all values $\alpha$ such that \ref{eq:CBB-CBA} holds, then together with any $\xi$, it contains $1-\xi$ and $\xi^k$ for any integer $k\geqslant 1$.
This implies that $\mathcal{A}$ is dense in $(0,1)$, and by continuity, we get $\mathcal{A}=(0,1)$.
That is, inequality \ref{eq:CBB-CBA} holds for any $\alpha\in (0,1)$.
This is equivalent to the point-on-side comparison for nonnegative curvature \cite[8.14]{AKP-2024}.
Similarly, one can show that if the inequality holds for some $\alpha>1$, then it holds for any $\alpha>1$,
and this is equivalent to the point-on-side comparison for nonpositive curvature \cite[9.14]{AKP-2024}.
\qeds

\section{Globalization}\label{par:globalization}

Let $K\subset W_n$ be a closed convex cone.
We say that a metric space $X$ satisfies \emph{local $K$-comparison} if any point $x\in X$ admits a neighborhood $U$ such that $K$-comparison holds for any $n$-point array in $U$.

If local $K$-comparison implies $K$-comparison for any length space, then we say that \emph{globalization holds} for $K$.

The following statement shows that if globalization holds for a nontrivial quadratic comparison, then it characterizes nonnegatively curved Alexandrov spaces;
so, Toponogov's theorem is the only nontrivial globalization theorem for quadratic conditions on 4-point arrays.

\begin{thm}{Theorem}\label{thm:globalization}
Suppose that globalization holds for a closed convex cone $K\z\subset W_4$.
Assume that the $K$-comparison is not trivial;
that is, on the one hand $\mathcal{M}_K$ does not include all length spaces and on the other hand $\mathcal{M}_K$ contains a space with at least two distinct points.
Then $\mathcal{M}_K$ consists of all Alexandrov spaces with nonnegative curvature.
\end{thm}

\begin{thm}{Lemma}\label{lem:globalization}
Under the assumptions of the theorem, $\mathcal{M}_K$ contains all Alexandrov spaces with nonnegative curvature.
\end{thm}

\parit{Proof.}
Since $K$-comparison is not trivial, $K\ne\{0\}$.
By \ref{prop:Associated form}, $\mathcal{M}_K$ contains the real line.
Therefore, local $K$-comparison holds
for any circle $r\cdot \mathbb{S}^1$ with $r>0$, and hence also for any product space $r\cdot \mathbb{S}^1\times\RR^3$.
It remains to apply \ref{prop:Four-point arrays}.
\qeds

In the following proof, we will use one statement from the next section.

\parit{Proof of the theorem.}
Let us denote by $K_0$ the cone in $W_4$ described by all inequalities of negative type $(3,1)$ and $(2,1)$.
By \ref{cor:Four-point arrays}, $K_0$ describes all metrics on 4-point arrays in Alexandrov spaces with nonnegative curvature.
By \ref{lem:globalization},  $K$~includes $K_0$.

Given a small $\delta>0$, consider all metrics with diameter at most $\delta$ on a 4-point array $\{x_1,x_2,x_3,x_4\}$ in an Alexandrov space with curvature at least $-1$.
Due to Wald's theorem (see \cite[Exercise 10.7]{AKP-2024}) we can assume that the $4$-point array
comes from a model $k$-plane (that is, a complete simply connected 2-dimensional Riemannian manifold of constant curvature $k$) with $k \geqslant -1$.
This implies that the array either (a) is isometric to a $4$-point subset of some sphere or (b) is bi-Lipschitz equivalent with constant $1 + C\cdot \delta^2$ to a $4$-point subset of a Euclidean plane, where $C > 0$ is some absolute constant.
Denote by $K_\delta\subset W_4$ the minimal closed convex cone that includes all associated forms of these metrics.
Note that $K_0\z\subset K_\delta$.
Furthermore, 
any element of $K_\delta$ is bi-Lipschitz equivalent with constant $1 + C \cdot\delta^2$ to some element of $K_0$;
therefore, $K_0=\bigcap_{\delta>0} K_\delta$.

Every Riemannian manifold admits a local curvature bound at each point.
In particular, after appropriate rescaling, a small neighborhood of any point of a Riemannian manifold has curvature at least $-1$.
Therefore, any compact Riemannian manifold satisfies the local $K_\delta$-comparison.

Suppose $K\supset K_\delta$ for some $\delta>0$.
Since globalization holds for $K$, the class $\mathcal{M}_K$ contains all compact Riemannian manifolds.
\textit{Every finite metric graph can be approximated by compact Riemannian manifolds.}
To prove it,
realize the graph in Euclidean space with smooth edges of the same length and take the boundary of an appropriate neighborhood.
(A much stronger result was proved by Vedrin Šahović in his thesis \cite{sahovic2009}.)
Any metric on $\{x_1,x_2,x_3,x_4\}$ admits a distance-preserving embedding into a metric graph, so $K$ contains a form near the associated form for any semimetric on $\{x_1,x_2,x_3,x_4\}$.
Since $K$ is closed, it contains forms associated with all metrics on the 4-point set;
so $K$ is defined only by the triangle inequalities, and the $K$-comparison is trivial.

From now on, we can assume that for any $\delta>0$, there is a form $\theta\in K_\delta\setminus K$.
By \ref{cor:squared-sides}, we can find a $(3,1)$-inequality that holds in $K$ with small error.
Namely, there is a $\lambda$-array $(\lambda_1,\lambda_2,\lambda_3,-1)$ such that $\lambda_i>0$ and
\[\begin{aligned}
\sum_{i,j}&\lambda_i\cdot\lambda_j\cdot|x_i-x_j|_X^2
\leqslant
\\
&\leqslant
10\cdot\delta^2\cdot \lambda_1\cdot\lambda_2\cdot\lambda_3\cdot (|x_1-x_2|_X^2+|x_2-x_3|_X^2+|x_3-x_1|_X^2)
\end{aligned}
\eqlbl{eq:+squares}\]
for any 4-point array $(x_1,x_2,x_3,x_4)$ with its form lying in $K$.

Let $\bm{\lambda}_\infty=(\lambda_1,\lambda_2,\lambda_3,-1)$ be a partial limit of these $\lambda$-arrays;
that is, for some sequence $\delta_n\to 0^+$, we can choose corresponding inequalities of the form \ref{eq:+squares} with $\lambda$-arrays $\bm{\lambda}_n$ such that $\bm{\lambda}_n\to \bm{\lambda}_\infty$ as $n\to \infty$.
We have to deal with three cases:
\begin{enumerate}
\item\label{in} $\lambda_1>0$, $\lambda_2>0$, and $\lambda_3>0$;
\item\label{side} $\lambda_i=0$ for one index $i$;
\item\label{vertex} $\lambda_i=0$ for two indices $i$.
\end{enumerate}

\parit{Case \ref{in}.}
Note that $(\lambda_1,\lambda_2,\lambda_3,-1)$ defines an inequality of negative type $(3,1)$.
This inequality holds for any form in $K$.
Therefore, \ref{prop:Alexandrov's comparison}\ref{prop:Alexandrov's comparison:(3,1)} finishes the proof.

\parit{Case \ref{side}.}
We can assume that $\lambda_3=0$, so
\[\bm{\lambda}_\infty=(\alpha,(1-\alpha),0,-1)\]
for some $0<\alpha<1$.
It defines an inequality of negative type $(2,1)$.
Note that
\[\bm{\lambda}_n=(\alpha_n\cdot(1-\beta_n),(1-\alpha_n)\cdot(1-\beta_n),\beta_n,-1),\]
for some sequences $\alpha_n\to\alpha$ and $\beta_n\to 0^+$ as $n\to\infty$.
A bit below we will show that repeating the proof of \ref{prop:Alexandrov's comparison}, we get
\[
\alpha\cdot|x_1-x_3|^2+(1-\alpha)\cdot|x_2-x_3|^2-\alpha\cdot (1-\alpha)\cdot|x_1-x_2|^2
\leqslant
|x_3-x_4|^2
\leqno{\text{\ref{eq:CBB-CBA}}'}
\]
if $x_4$ lies on $[x_1x_2]$ and divides it in the ratio $(1-\alpha):\alpha$.
After that it remains to follow the end of the proof of \ref{prop:Alexandrov's comparison}.

To prove \ref{eq:CBB-CBA}$'$, suppose $x_4$ lies on $[x_1x_2]$ and divides it in the ratio $(1-\alpha_n):\alpha_n$.
Then
\[\alpha_n\cdot (1-\alpha_n)\cdot|x_1-x_2|^2-\alpha_n\cdot |x_1-x_4|^2-(1-\alpha_n)\cdot |x_2-x_4|^2=0;\leqno{\text{\ref{eq:trig-inq}}'}\]
this will be used instead of \ref{eq:trig-inq} (page \pageref{eq:trig-inq}).
Set
\begin{align*}
S&=|x_1-x_2|_X^2+|x_2-x_3|_X^2+|x_3-x_1|_X^2,
\\
E_n&=10\cdot \alpha_n\cdot (1-\alpha_n)\cdot(1-\beta_n)^2\cdot \beta_n\cdot S.
\end{align*}
If we substitute $\alpha=\alpha_n$ and $\gamma=\beta_n$ in \ref{eq:lambda-inq}, then instead of zero in the right-hand side, we get $\delta_n^2\cdot E_n$:
\[
\begin{aligned}
\alpha_n\cdot (1-\alpha_n)\cdot(1-\beta_n)^2\cdot|x_1-x_2|^2 - \beta_n\cdot |x_3-x_4|^2 &+
\\
+(1-\alpha_n)\cdot(1-\beta_n)\cdot\beta_n\cdot|x_2-x_3|^2-\alpha_n\cdot (1-\beta_n)\cdot |x_1-x_4|^2&+
\\
+\alpha_n\cdot(1-\beta_n)\cdot\beta_n\cdot|x_1-x_3|^2-(1-\alpha_n)\cdot(1-\beta_n)\cdot |x_2-x_4|^2&\leqslant \delta_n^2\cdot E_n
\end{aligned}
\leqno{\text{\ref{eq:lambda-inq}}'}
\]

The inequality $\tfrac1{\beta_n}\cdot$\ref{eq:lambda-inq}$'\,-\,\tfrac{1-\beta_n}{\beta_n}\cdot$\ref{eq:trig-inq}$'$ implies the following
\[
\begin{aligned}
\alpha_n\cdot|x_1-x_3|^2+(1-\alpha_n)\cdot|x_2-x_3|^2&-
\alpha_n\cdot (1-\alpha_n)\cdot|x_1-x_2|^2 \leqslant
\\
&\leqslant |x_3-x_4|^2 +  \tfrac{\delta_n^2}{\beta_n}\cdot E_n+\beta_n\cdot S.
\end{aligned}
\]
It remains to observe that $\tfrac{\delta_n^2}{\beta_n}\cdot E_n+\beta_n\cdot S\to0$ as $n\to\infty$.

\parit{Case \ref{vertex}.}
We may assume that $\lambda_1=\lambda_2=0$, so $\bm{\lambda}_\infty=(0,0,1,-1)$ and
\[\bm{\lambda}_n=(\alpha_n\cdot(1-\beta_n),(1-\alpha_n)\cdot(1-\beta_n),\beta_n,-1),\]
where $0<\alpha_n<1$, $0<\beta_n<1$, and $\beta_n\to 1^-$ as $n\to\infty$.

The following argument pushes $\bm{\lambda}_\infty$ out of the corner; that is, it produces a new sequence of inequalities~\ref{eq:L<En-hats} with new arrays ${\bm{\lambda}}'_n$ that fall under case~\ref{in} or~\ref{side}.

We start with a sequence of inequalities
\[
L_n(x_1,x_2,x_3,x_4)\leqslant \delta_n^2\cdot E_n(x_1,x_2,x_3),
\eqlbl{eq:L<En}
\]
where $\delta_n\to0$ as $n\to\infty$,
$L_n$ is the left-hand side of~\ref{eq:lambda} for $\bm{\lambda}_n$, and
\[
E_n=10\cdot\alpha_n\cdot(1-\alpha_n)\cdot\beta_n\cdot(1-\beta_n)^2
\cdot (|x_1-x_2|_X^2+|x_2-x_3|_X^2+|x_3-x_1|_X^2).
\]

Choose the largest $k = k(n)$ such that $ \beta_n^{k}\geqslant \tfrac{99}{100}$.
Note that $k(n)$ is defined for all sufficiently large~$n$;
moreover, $\beta_n^k\to \tfrac{99}{100}$ as $n\to\infty$.
As in the proof of Claim~\ref{clm:1=>2}, let $L_n'$ be the left-hand side in~\ref{eq:lambda} for the array
\[{\bm{\lambda}}'_n=(\alpha_n\cdot(1-\beta_n^{k}),(1-\alpha_n)\cdot(1-\beta_n^{k}),\beta_n^{k},-1),\]
and define
\[E_n'=10\cdot\alpha_n\cdot(1-\alpha_n)\cdot\beta_n^{k}\cdot(1-\beta_n^{k})^2 \cdot (|x_1-x_2|_X^2+|x_2-x_3|_X^2+|x_3-x_1|_X^2).\]

Let us prove the following inequality
\[
L_n'\leqslant 1000\cdot \delta_n^2\cdot E_n'.
\eqlbl{eq:L<En-hats}
\]
Note that the new array ${\bm{\lambda}}'_n$ falls under case~\ref{in} or~\ref{side}.
Therefore, once \ref{eq:L<En-hats} is proved, case \ref{vertex} will be settled.

Assume the contrary; choose an array $x_1,x_2,x_3,x_{3+k}$ for which \ref{eq:L<En-hats} does not hold.
Then
\[
|x_3-x_{3+k}|\leqslant \max\{\,|x_3-x_1|,|x_3-x_2|\,\},
\]
since otherwise $L_n'(x_1,x_2,x_3,x_{3+k})\leqslant 0$, which would imply~\ref{eq:L<En-hats}.

We can assume that $X$ is geodesic; otherwise, pass to its ultracompletion.
By \ref{eq:L'=L+L},
$L_n'(x_1,x_2,x_3,x_{3+k})$ is a linear combination
of $L_n(x_1,x_2,x_{2+i},x_{3+i})$ with positive coefficients;
here $x_3,\ldots,x_{3+k}$ are points on the geodesic $[x_3,x_{3+k}]$.
Note that the coefficient in front of each $L_n(x_1,x_2,x_{2+i},x_{3+i})$ does not exceed $k$.
Therefore,
\[L_n'(x_1,x_2,x_3,x_{3+k})\leqslant k\cdot \delta_n^2\cdot \sum_{i=1}^{k}E_n(x_1,x_2,x_{2+i}).
\]

Since $|x_3-x_{3+k}|\leqslant \max\{\,|x_3-x_1|,|x_3-x_2|\,\}$, we have
\begin{align*}
|x_1-x_{2+i}|&\leqslant 2\cdot \max\{\,|x_3-x_1|,|x_3-x_2|\,\},\\
|x_2-x_{2+i}|&\leqslant 2\cdot \max\{\,|x_3-x_1|,|x_3-x_2|\,\},
\end{align*}
for~$i\geqslant1$.
It follows that
\[
E_n(x_1,x_2,x_{2+i})\leqslant 100\cdot E_n(x_1,x_2,x_3)
\]
for all~$i$; hence
\[L_n'(x_1,x_2,x_3,x_{3+k})\leqslant 100\cdot k^2\cdot \delta_n^2\cdot E_n(x_1,x_2,x_3).
\]
Furthermore, since $1>\beta_n\geqslant\beta_n^{k}\geqslant \tfrac{99}{100}$,
we have $1-\beta_n^{k}\geqslant \tfrac k2\cdot(1-\beta_n)$, and hence
\[
k^2\cdot E_n(x_1,x_2,x_3)\leqslant 10\cdot E_n'(x_1,x_2,x_3)
\]
for all large $n$.
It follows that \ref{eq:L<En-hats} holds for $x_1,x_2,x_3,x_{3+k}$ --- a contradiction.
\qeds

\section{Auxiliary statements}\label{Auxiliary statements}

\begin{thm}{Lemma}\label{lem:area-bound}
Let $(\lambda_1,\lambda_2,\lambda_3,-1)$ be the $\lambda$-array of an inequality of type $(3,1)$ and
let $(x_1$, $x_2$, $x_3$, $x_4)$ be a 4-point array in the hyperbolic space.

Then
\[\sum_{i,j=1}^4\lambda_i\cdot\lambda_j\cdot|x_i-x_j|_X^2
\leqslant
(24+6\cdot\delta^2)\cdot\lambda_1\cdot\lambda_2\cdot\lambda_3\cdot\tilde a^2,\]
where $\delta=\max\{|x_1-x_2|,|x_1-x_3|,|x_2-x_3|\}$
and
$\tilde a$ is the area of the model triangle $\tilde\triangle(x_1x_2x_3)_{\EE^2}$.
\end{thm}

Note that Heron's formula implies
\begin{align*}
16\cdot \tilde a^2
\quad=\quad &(|x_1-x_2|^2+|x_2-x_3|^2+|x_3-x_1|^2)^2-
\\
-2\cdot &(|x_1-x_2|^4+|x_2-x_3|^4+|x_3-x_1|^4).
\end{align*}
In particular, $\tilde a^2$ is a quadratic form on $W_4$.

\parit{Proof.}
By Kirszbraun's theorem (see \cite{lang-schroeder,AKP-2011} or \cite[Chapter 10]{AKP-2024}), we may assume that $x_1$, $x_2$, $x_3$ and $x_4$ lie in the hyperbolic plane; moreover, $x_4$ belongs to the solid hyperbolic triangle with vertices $x_1$, $x_2$ and $x_3$.
Denote by $a$ the area of this triangle.
Again, by Kirszbraun's theorem,
\[a\leqslant \tilde a.\]

Since the hyperbolic plane has curvature $-1$, we have
\[\pi-\angk{x_1}{x_2}{x_3}_{\HH^2}-\angk{x_2}{x_3}{x_1}_{\HH^2}-\angk{x_3}{x_1}{x_2}_{\HH^2}=a.\]
These relations and the comparison imply that
\begin{align*}
0&\leqslant \angk{x_1}{x_2}{x_3}_{\EE^2}-\angk{x_1}{x_2}{x_3}_{\HH^2}\leqslant \tilde a.
\intertext{Since $x_4$ lies in the solid hyperbolic triangle with vertices $x_1$, $x_2$ and $x_3$, we also have $\angk{x_1}{x_2}{x_4}_{\HH^2}+\angk{x_1}{x_4}{x_3}_{\HH^2}=\angk{x_1}{x_2}{x_3}_{\HH^2}$.
Hence, the comparison also implies}
0&\leqslant
\angk{x_1}{x_2}{x_4}_{\EE^2}+\angk{x_1}{x_4}{x_3}_{\EE^2}-\angk{x_1}{x_2}{x_3}_{\HH^2}\leqslant \tilde a.
\end{align*}

Set $\phi=\angk{x_1}{x_2}{x_3}_{\EE^2}$ and $\psi=\angk{x_1}{x_2}{x_4}_{\EE^2}+\angk{x_1}{x_4}{x_3}_{\EE^2}$.
From the above, we get $\phi\leqslant \psi+\tilde a$.
By the law of cosines,
\[|x_2-x_3|^2=|x_1-x_2|^2+ |x_1-x_3|^2-2|x_1-x_2|\cdot|x_1-x_3|\cdot\cos\phi.\]
Redefining the distance $|x_2-x_3|$ via the law of cosines with angle $\psi$,
\[|x_2-x_3|^2\mathrel{:=}|x_1-x_2|^2+ |x_1-x_3|^2-2|x_1-x_2|\cdot|x_1-x_3|\cdot\cos\psi,\]
yields a Euclidean quadruple.
That is, decreasing $|x_2-x_3|^2$ by
\[s=2\cdot |x_1-x_2|\cdot|x_1-x_3|\cdot(\cos\psi-\cos\phi)\]
makes the quadruple Euclidean.
Since
\[\cos\psi-\cos\phi=\sin\phi\cdot\sin(\phi-\psi)-2\cdot\cos\phi\cdot (\sin\tfrac{\phi-\psi}2)^2,\] we get
\begin{align*}
s
&\leqslant\
2\cdot|x_1-x_2|\cdot|x_1-x_3|\cdot (\tilde a\cdot\sin\phi+\tfrac12\cdot\tilde a^2)
\leqslant
\\
&\leqslant\ (4+\delta^2)\cdot \tilde a^2.
\end{align*}

It follows that
\[\sum_{i,j}\lambda_i\cdot\lambda_j\cdot|x_i-x_j|_X^2\leqslant 2\cdot\lambda_2\cdot\lambda_3\cdot(4+\delta^2)\cdot \tilde a^2.\]
We may assume that $\lambda_1\geqslant \lambda_2\geqslant \lambda_3>0$.
Since $\lambda_1+ \lambda_2+ \lambda_3=1$, we have $\lambda_1\geqslant \tfrac13$, and the statement follows.
\qeds

Recall that $K_0$ denotes the cone in $W_4$ described by all inequalities of negative type $(3,1)$ and $(2,1)$.

\begin{thm}{Corollary}\label{cor:squared-sides}
Let $K\subset W_4$ be a closed convex cone.
Suppose that $K\supset K_0$, but $K$-comparison does not hold for a 4-point array of diameter at most $\delta$ in the hyperbolic space.
If $\delta$ is sufficiently small,
then there is a $\lambda$-array $(\lambda_1,\lambda_2,\lambda_3,-1)$ of type $(3,1)$
such that the inequality
\[\sum_{i,j=1}^4\lambda_i\cdot\lambda_j\cdot|x_i-x_j|_X^2
\leqslant
10\cdot\delta^2\cdot \lambda_1\cdot\lambda_2\cdot\lambda_3\cdot (|x_1-x_2|_X^2+|x_2-x_3|_X^2+|x_3-x_1|_X^2)\]
holds for any 4-point array whose associated form lies in $K$.
\end{thm}

This inequality should be interpreted as a $(3,1)$-inequality with error
\[10\cdot\delta^2\cdot \lambda_1\cdot\lambda_2\cdot\lambda_3\cdot (|x_1-x_2|_X^2+|x_2-x_3|_X^2+|x_3-x_1|_X^2).\]

\parit{Proof.}
Let $\theta$ be the associated form of the 4-point array in the hyperbolic space.
Since $K$ is convex, we can choose a quadratic inequality that does not hold for $\theta$, but holds for any form in $K$.
Our inequality can be written as $\langle \omega,\rho_{\bm{x}} \rangle\geqslant 0$, where $\rho_{\bm{x}}\in W_4$ is an associated form of some 4-point array $(x_1,x_2,x_3,x_4)$ and $\omega$ is a fixed unit form in $W_4$.

\begin{wrapfigure}{o}{30mm}
\centering
\vskip-4mm
\includegraphics{mppics/pic-40}
\vskip-0mm
\end{wrapfigure}

Let $N\subset \SSS^2\subset V_4$ be the set of unit vectors $v$ for which $\theta(v)<0$.
The set $N$ does not intersect the four equators $e_1$, $e_2$, $e_3$, $e_4$ parallel to the facets of $\triangle$.
Moreover, since $K\supset K_0$, the set $N$ (up to sign) lies in one of the triangles, say $T$;
we can assume that equators $e_1$, $e_2$, $e_3$ are extensions of sides of $T$, and $e_4$ is the remaining equator.

Since $K\supset K_0$, we have $\omega(v)\geqslant 0$ for any $v\in V_4$.
Let $\sigma$ be a unit 1-form on $V_4$ that vanishes on $e_4$.
By the spectral theorem, $\omega$ and $\sigma^2$ can be diagonalized in a common basis (which does not have to be orthogonal).
Therefore,
\[\langle \omega,\rho_{\bm{x}}\rangle=\rho_{\bm{x}}(u)+\rho_{\bm{x}}(v)+\rho_{\bm{x}}(w)\eqlbl{eq:v_123}\]
for fixed vectors $u,v,w\in V_4$ such that $v$ and $w$ lie in the plane spanned by $e_4$;
the latter condition follows since $\sigma^2$ is also diagonalized.

The inequalities $\rho_{\bm{x}}(v)\geqslant 0$ and
$\rho_{\bm{x}}(w)\geqslant 0$ are of type $(2,1)$;
they follow from the triangle inequality, so we always have that $\rho_{\bm{x}}(v)\geqslant 0$ and
$\rho_{\bm{x}}(w)\geqslant 0$.
Moreover, the values $\rho_{\bm{x}}(v)$ and $\rho_{\bm{x}}(w)$ depend only on the sides of the triangle $[x_1x_2x_3]$.

Since
\[
\begin{aligned}
\theta(v)&\geqslant0,
&
\theta(w)&\geqslant0,
&
\text{and}&&
\langle\omega,\theta\rangle=\theta(u)+\theta(v)+\theta(w)&<0,
\end{aligned}
\eqlbl{eq:uvw}
\]
we have that $u$ points in $N$.
Therefore, $\rho_{\bm{x}}(u)\geqslant0$ is an inequality of negative type $(3,1)$.

Let us rescale $u$, $v$, and $w$ so that the inequality $\rho_{\bm{x}}(u)\geqslant 0$ (which as we know has negative type $(3,1)$)
has $\lambda$-array $(\lambda_1,\lambda_2,\lambda_3,-1)$.
After this rescaling, the left-hand side in the required inequality can be written as $-\rho_{\bm{x}}(u)$.
Now we have
\[\sum_{i,j=1}^4\lambda_i\cdot\lambda_j\cdot|x_i-x_j|_X^2=-\rho_{\bm{x}}(u)\leqslant \rho_{\bm{x}}(v)+\rho_{\bm{x}}(w),\]
and it remains to estimate $\rho_{\bm{x}}(v)+\rho_{\bm{x}}(w)$.

Since the inequalities $\rho_{\bm{x}}(v)\geqslant 0$ and $\rho_{\bm{x}}(w)\geqslant 0$ have type $(2,1)$, we can choose $\varepsilon_v\geqslant 0$ and $\varepsilon_w\geqslant 0$ such that
\[\rho_{\bm{x}}(v)=\varepsilon_v\cdot \tilde d_v^2
\quad\text{and}\quad
 \rho_{\bm{x}}(w)=\varepsilon_w\cdot \tilde d_w^2,
\]
where
$\tilde d_v$ and $\tilde d_w$ are distances from a vertex of the model triangle $\tilde\triangle(x_1x_2x_3)_{\EE^2}$ to a point on the opposite side that divides it in a certain ratio (the choice of vertex and ratio depends on $v$ and $w$, respectively).
Note that $\tilde d_v$ and $\tilde d_w$ do not exceed the largest side of $\tilde\triangle(x_1x_2x_3)_{\EE^2}$; in particular
\begin{align*}
\tilde d_v^2&\leqslant |x_1-x_2|_X^2+|x_2-x_3|_X^2+|x_3-x_1|_X^2,
\\
\tilde d_w^2&\leqslant |x_1-x_2|_X^2+|x_2-x_3|_X^2+|x_3-x_1|_X^2.
\end{align*}
Therefore, it is sufficient to find an appropriate upper bound on $\varepsilon_v+\varepsilon_w$.

Now suppose that the array $\bm{x}=(x_1,x_2,x_3,x_4)$ is equipped with a metric such that its associated form is $\theta$.
Since $\delta$ is small, \ref{eq:uvw} and \ref{lem:area-bound} imply that
\[\varepsilon_v\cdot \tilde d_v^2+\varepsilon_w\cdot \tilde d_w^2=\theta(v)+\theta(w)
<
-\theta(u)
\leqslant
25\cdot\lambda_1\cdot\lambda_2\cdot\lambda_3\cdot\tilde a^2.\]
Since the diameter of $\tilde\triangle(x_1x_2x_3)$ is at most $\delta$, we also have
$\tilde a\z\leqslant \delta\cdot \tilde d_v/2$ and
$\tilde a\z\leqslant \delta\cdot \tilde d_w/2$.
Hence
\[\varepsilon_v+\varepsilon_w\leqslant \tfrac{25}4\cdot \delta^2\cdot\lambda_1\cdot\lambda_2\cdot\lambda_3,\]
which implies the required inequality.
\qeds

\section{Remarks}

The requirement in \ref{thm:globalization} that the cone $K$ be closed is necessary;
without it, the globalization holds for the cone described by the inequality
\[|x_1-x_2|^2+|x_1-x_3|^2+|x_2-x_3|^2<4\cdot(|x_1-x_4|^2+|x_2-x_4|^2+|x_3-x_4|^2),\]
assuming that the right-hand side is positive.
To verify this statement, note that the inequality forbids tripods.

It would be super-nice to find a quadratic inequality for $n$-point arrays with globalization that is not implied by the standard globalization theorem.
Some candidates can be found in our earlier paper \cite{lebedeva-petrunin-zolotov}.
Also, it might be possible to prove a version of our globalization theorem
for simply connected spaces, which includes the Cartan--Hadamard theorem.

One may consider arbitrary closed conditions on $n$-point arrays.
It would be even more interesting to understand the globalization phenomenon in such a general setting.
One family of examples is provided by the conditions of curvature bounded below by $\kappa\in\RR$ in the sense of Alexandrov; however, globalization also holds for other conditions.
For example, it holds for the condition determined by the set of semimetrics on four-point arrays that admit an isometric embedding into $r\cdot\SSS^1$ for some $r>0$.

\paragraph{Acknowledgments.}
Nina Lebedeva was supported by the Ministry of Science and Higher Education of the Russian Federation, agreement no. 075-15-2025-343;
Anton Petrunin was supported by NSF grant DMS-2005279.